\documentclass[12pt]{article}

 \usepackage{a4}
 \usepackage{amsmath}
 \usepackage{amssymb}
 \usepackage{graphicx}
\newtheorem{thm}{Theorem}[section]

 \newtheorem{rem}[thm]{Remark}

 \renewcommand{\a}{\alpha}

\title{On an iteration leading to a $q$-analogue of the Digamma function}
 \author{Christian Berg\footnote{Corresponding author}\; and 
Helle Bjerg Petersen} 
 \date{\today}

\begin{document}

 \maketitle 

\begin{abstract} We show that the $q$-Digamma function $\psi_q$ for
  $0<q<1$ appears in an
  iteration studied by Berg and Dur{\'a}n. In addition we determine
  the probability measure $\nu_q$ with moments $1/\sum_{k=1}^{n+1}
  (1-q)/(1-q^k)$, which are  $q$-analogues of the reciprocals of the
  harmonic numbers.
\end{abstract}


    
\section{Introduction}

For a measure $\mu$ on the unit interval $[0,1]$ we consider
its {\it Bernstein transform}
\begin{equation}\label{eq:Bern}
\mathcal B(\mu)(z)=\int _0^1\frac{1-t^z}{1-t}d\mu (t),\quad  \Re z> 0,
\end{equation}
as well as its {\it Mellin transform}
\begin{equation}\label{eq:Mellin}
\mathcal M(\mu)(z)=\int _0^1t^zd\mu (t),\quad  \Re z> 0.
\end{equation}
These functions are clearly holomorphic in
the right half-plane $\Re z>0$.

The two integral transformations are combined in the following
theorem from \cite{B:D1} about Hausdorff moment sequences, i.e.,
sequences $(a_n)_{n\ge 0}$ of the form
\begin{equation}\label{eq:Haus}
a_n=\int_0^1 t^n\,d\mu(t),
\end{equation}
for a positive measure $\mu$ on the unit interval.

\begin{thm} \label{thm:HtoH}
Let $(a_n)_{n\ge 0}$ be a Hausdorff moment sequence as in
\eqref{eq:Haus} with $\mu\neq 0$. Then the sequence
$(T(a_n))_{n\ge 0}$ defined by $T(a_n)_n=1/(a_0+\ldots +a_n)$ is again a
Hausdorff moment sequence, and its associated measure $\widehat{T}(\mu
)$ has the properties $\widehat{T}(\mu)(\{ 0\})=0$ and
\begin{equation}\label{eq:assomeas}
\mathcal B(\mu)(z+1)\mathcal M(\widehat{T}(\mu))(z)=1\quad
\mbox{for}\quad \Re z> 0.
\end{equation}
\end{thm}
This means that the measure $\widehat{T}(\mu)$ is determined as
the inverse Mellin transform of the function $1/\mathcal
B(\mu)(z+1)$.

It follows by Theorem \ref{thm:HtoH} that $T$ maps the set of
normalized Hausdorff moment sequences (i.e., $a_0=1$) into itself.
By Tychonoff's extension of Brouwer's fixed
point theorem, $T$ has a fixed point $(m_n)$. Furthermore, it is clear that a
fixed point $(m_n)$ is uniquely determined by the equations
\begin{equation}\label{eq:fixed}
(1+m_1+\ldots +m_n)m_n=1,\quad n\ge 1.
\end{equation}                                    
Therefore
\begin{equation}\label{eq:fixrec}
m_{n+1}^2+\frac{m_{n+1}}{m_n}-1=0,
\end{equation}
giving
$$
m_1=\frac{-1+\sqrt 5}{2},\quad m_2=\frac{\sqrt{22+2\sqrt 5}-\sqrt
5-1}{4}, \ldots\,.
$$  

Similarly, $\widehat{T}$ maps the set $M_+^1([0,1])$ of probability
measures on $[0,1]$ into itself. It has a uniquely determined fixed
point $\omega$ and 
\begin{equation}\label{eq:fix}
m_n=\int_0^1 t^n\,d\omega(t),\quad n=0,1,\ldots.
\end{equation}

 Berg and Dur{\'a}n studied this fixed point in \cite{B:D2},\cite{B:D3}, and it was
proved that the Bernstein transform $f=\mathcal B(\omega)$ is
meromorphic in the whole complex plane and characterized by a
functional equation  and a log-convexity property in ana\-lo\-gy with
Bohr-Mollerup's characterization of the Gamma 
function. Let us  also mention that $\omega$ has an increasing and
convex  density with respect to Lebesgue measure $m$ on the unit interval. 

An important step in the proof is to establish that $\omega$ is an
attractive fixed point so that in particular the iterates
$\widehat{T}^{\circ n}(\delta_1)$ converge weakly to $\omega$. Here and
in the following $\delta_a$ denotes the Dirac measure with mass 1
concentrated in $a\in \mathbb R$.

It is easy to see that $\widehat{T}(\delta_1)=m$, because 
$$
T(1,1,\ldots)_n=\frac{1}{n+1}=\int_0^1 t^n\,dt.
$$

It is well-known that the Bernstein transform of  Lebesgue measure
$m$ on $[0,1]$ is related to the Digamma function $\psi$, i.e., the
logarithmic derivative of the Gamma function, since
\begin{equation}\label{eq:digamma}
\int _0^1\frac{1-t^z}{1-t}dt=\psi(z+1)+\gamma=\sum_{n=1}^\infty
\frac{z}{n(n+z)},\quad \Re z>0,
\end{equation} 
cf. \cite[8.36]{Gr:Ry}. Here $\gamma=-\psi(1)$ is Euler's constant.

 Therefore
$\nu_1:=\widehat{T}(m)=\widehat{T}^{\circ 2}(\delta_1)$ is determined
by
$$
\mathcal M(\nu_1)(z)=\frac{1}{\mathcal B(m)(z+1)}=\frac{1}{\psi(z+2)+\gamma}.
$$

The measure $\nu_1=\widehat{T}(m)$ is given explicitly in \cite{B:D1} as
\begin{equation}\label{eq:nu1}
\nu_1=(\sum_{n=0}^\infty \a_nt^{\xi_n})\,dt,
\end{equation}
where $\xi_0=0$, $\xi_n\in (n,n+1),n=1,2,\ldots$ is the solution
to $\psi(1-\xi_n)=-\gamma$ and $\a_n=1/\psi'(1-\xi_n)$. The moments of
the measure $\nu_1$ are the reciprocals of the harmonic numbers, i.e.,
\begin{equation}\label{eq:harmonic}
\int_0^1 t^n\;d\nu_1(t)=\frac{1}{\mathcal H_{n+1}}=\left(\sum_{k=1}^{n+1}\frac{1}{k}\right)^{-1}.
\end{equation}

The purpose of this paper is to study the first elements of the
sequence $\widehat{T}^{\circ n}(\delta_q)$, where $0<q<1$ is
fixed. The reason for excluding $q=0$ is that
$\widehat{T}(\delta_0)=\delta_1$. Since $\omega$ is an attractive fixed point, we know that the
sequence converges weakly to $\omega$. 

The first step in the iteration is easy: 
\begin{equation}\label{eq:deltaq}
\widehat{T}(\delta_q)=(1-q)\sum_{k=0}^\infty q^k\delta_{q^k},
\end{equation}
because
\begin{equation}\label{eq:deltaq1}
\int_0^1 t^z
\,d\,\widehat{T}(\delta_q)(t)=\frac{1-q}{1-q^{z+1}}=(1-q)\sum_{k=0}^\infty
q^k q^{kz}.
\end{equation}

This shows that $\widehat{T}(\delta_q)$ is the Jackson $d_qt$-measure on
$[0,1]$ used in the theory of $q$-integrals, cf. \cite{G:R}. It is a
$q$-analogue of Lebesgue measure in the sense that $d_qt\to m$ weakly
for $q\to 1$. 

It is therefore to be expected that $\nu_q:=\widehat{T}(d_qt)=\widehat{T}^{\circ 2}(\delta_q)$ is a
$q$-analogue of the measure $\nu_1$, and we are going to determine
$\nu_q$ as closely as possible. We have
\begin{equation}\label{eq:Melnuq}
\mathcal M(\nu_q)(z)=\frac{1}{f_q(z+1)},
\end{equation}
where $f_q$ is defined as the Bernstein transform of $d_qt$:
\begin{equation}\label{eq:Bdqt}
f_q(z)=\int_0^1 \frac{1-t^z}{1-t}\,d_qt=(1-q)\left(z+\sum_{k=1}^\infty
q^k \frac{1-q^{kz}}{1-q^k}\right).
\end{equation} 
This formula is a $q$-analogue of \eqref{eq:digamma}.
The moments of $\nu_q$ are $q$-analogues of \eqref{eq:harmonic} 
\begin{equation}\label{eq:q-harmonic}
\int_0^1 t^n\,d\nu_q(t)=\left(\sum_{k=0}^{n}\frac{1-q}{1-q^{k+1}}\right)^{-1}.
\end{equation}

Our main result is the following (note that the Haar measure on the
multiplicative group $]0,\infty[$ is $dt/t$):

\begin{thm}\label{thm:main} The measure $\nu_q$ has a continuous
  density $\widetilde{\nu_q}(t)$ with  respect to $dt/t$ on
  $]0,1]$. It is $C^\infty$ on each of the open intervals
  $]q^{k+1},q^k[,k=0,1,\ldots$ with jump of the derivative of size
  $q^k/(1-q^k)(1-q)$ at the point $q^k,k=1,2,\ldots$. Furthermore,
  $\lim_{t\to 0}\widetilde{\nu_q}(t)=0$.
\end{thm}

\begin{rem} {\rm It follows that the behaviour of
    $\widetilde{\nu_q}(t)$ is oscillatory, and therefore quite
    different from that of $\lim_{q\to 1}\widetilde{\nu_q}(t)$, which
is increasing and convex. In fact, it follows from \eqref{eq:nu1} that
$$
\lim_{q\to 1}\widetilde{\nu_q}(t)=\sum_{n=0}^\infty \alpha_n
t^{\xi_n+1},\quad 0<t\le 1.
$$
See Figure 1 and 2 which shows the graph of $\widetilde{\nu_q}(e^{-t})$ for
$q=0.5$ and $q=0.9$.
}
\end{rem}
\section{Proofs}

Jackson's $q$-analogue of the Gamma function is defined as
$$
\Gamma_q(z)=\frac{(q;q)_\infty}{(q^z;q)_\infty}(1-q)^{1-z},
$$
cf. \cite{G:R}, and its logarithmic derivative
\begin{equation}\label{eq:psiq}
\psi_q(z)=\frac{d}{dz}\log\Gamma_q(z)=-\log(1-q)+\log q\sum_{k=0}^\infty\frac{q^{k+z}}{1-q^{k+z}}
\end{equation}
has been proposed in \cite{K:S} as a $q$-analogue of the Digamma
function $\psi$. See also the recent paper \cite{M:S}.
We define the $q$-analogue of Euler's constant as
\begin{equation}\label{eq:qgamma}
\gamma_q=-\psi_q(1)=\log(1-q)-\log q\sum_{k=1}^\infty\frac{q^k}{1-q^k}.
\end{equation} 
The  Bernstein transform $f_q$ of $d_qt$ is given in \eqref{eq:Bdqt}, hence
\begin{eqnarray*}
\frac{f_q(z)}{1-q} &=& z+\sum_{k=1}^\infty q^k\sum_{n=0}^\infty q^{kn}(1-q^{kz})\\ 
       &=& z+\sum_{n=0}^\infty\left(\sum_{k=1}^\infty
           (q^{k(n+1)}-q^{k(n+1+z)})\right)\\
       &=& z+ \frac{1}{\log(1/q)}\left(\gamma_q+\psi_q(z+1)\right),
\end{eqnarray*}
which shows the close relationship with the $q$-Digamma function.
We will be using another expression for $f_q(z)/(1-q)$ derived from
\eqref{eq:Bdqt}, namely
\begin{equation}\label{eq:Bdqt1} 
\frac{f_q(z)}{1-q}=z+c_q-\sum_{k=1}^\infty \frac{q^k}{1-q^k}q^{kz},
\end{equation}
with
\begin{equation}\label{eq:cq}
c_q=\sum_{k=1}^\infty\frac{q^k}{1-q^k}.
\end{equation}
Clearly, $q/(1-q)<c_q<q/(1-q)^2$ for $0<q<1$ and $q\mapsto c_q$ is a
stricly increasing  map of $]0,1[$ onto $]0,\infty[$. We mention two
other expressions 
$$
c_q=\sum_{n=1}^\infty d(n)q^n=\sum_{n=1}^\infty\left(1-(q^n;q)_\infty\right),
$$
where $d(n)$ is the number of divisors in $n$, see \cite[p. 14]{F}.

In order to replace  the Mellin transformation by the Laplace
transformation we introduce the probability measure $\tau_q$ on
$[0,\infty[$ which has $\nu_q$ as image measure under $t\to e^{-t}$,
hence
$$
\mathcal L(\tau_q)(z)=\int_0^\infty
e^{-tz}\,d\tau_q(t)=\frac{1}{f_q(z+1)}.
$$
The analogue of Theorem~\ref{thm:main} about the measure $\tau_q$ is
given in the next theorem, which we shall prove first.

\begin{thm}\label{thm:tauq} The measure $\tau_q$ has a continuous
  density also denoted $\tau_q$ with respect to Lebesgue measure on
  $[0,\infty[$. It is $C^\infty$ in each of the open intervals
  $]n\log(1/q),(n+1)\log(1/q)[,n=0,1,\ldots$ with jump of the
  derivative of size
\begin{equation}\label{eq:jump}
J_n=\frac{q^{2n}}{(1-q^n)(1-q)}
\end{equation}
at the point
  $n\log(1/q),n=1,2,\ldots$. Furthermore, $\lim_{t\to\infty}\tau_q(t)=0$.
\end{thm} 

{\it Proof of Theorem~\ref{thm:tauq}}. Introducing the discrete measure
$$
\mu=\sum_{k=1}^\infty \frac{q^{2k}}{1-q^k}\delta_{k\log(1/q)}
$$
of finite total mass 
\begin{equation}\label{eq:total}
||\mu||_1=c_q-q/(1-q)<c_q,
\end{equation}
we can write
$$
\frac{f_q(z+1)}{1-q}=1+c_q+z-\mathcal L(\mu)(z),
$$
hence
\begin{equation}\label{eq:factor}
\frac{1-q}{f_q(z+1)}=\left((1+c_q+z)(1-\frac{\mathcal
    L(\mu)(z)}{1+c_q+z})\right)^{-1}=\sum_{n=0}^\infty\frac{\left(\mathcal
  L(\mu)(z)\right)^n}{\left(1+c_q+z\right)^{n+1}}.
\end{equation}
Let $\rho_q$ denote the following exponential density restricted to
the positive half-line 
$$
\rho_q(t)=\exp(-(1+c_q)t)Y(t),
$$
where $Y$ is the usual Heaviside function equal to 1 for $t\ge 0$ and
equal to zero for $t<0$. Its Laplace transform is given as
$$
\int_0^\infty e^{-tz}\rho_q(t)\,dt=(1+c_q+z)^{-1},
$$
but this shows that \eqref{eq:factor} is equivalent to the  following convolution equation
\begin{equation}\label{eq:conv}
(1-q)\tau_q=\rho_q\ast\sum_{n=0}^\infty(\mu\ast \rho_q)^{\ast
  n}=\sum_{n=0}^\infty \rho_q^{\ast(n+1)}\ast\mu^{\ast n}.
\end{equation}

This equation expresses a factorization of $(1-q)\tau_q$  as the convolution of
the exponential density $\rho_q$ and an elementary kernel
$\sum_0^\infty N^{\ast n}$ with $N=\mu\ast\rho_q$. For
information about the basic notion of elementary kernels in potential theory,
see \cite[p.100]{B:F}.
All three measures in question $\tau_q,\rho_q$ and $\sum_0^\infty
(\mu\ast\rho_q)^{\ast n}$ are potential kernels on $\mathbb R$ in the
sense of \cite{B:F}.

The measure $\mu^{\ast n},n\ge 1$ is a discrete measure concentrated
in the points  $k\log(1/q),k=n,n+1,\ldots$. The convolution powers of
$\rho_q$ are Gamma densities
$$
\rho_q^{\ast (n+1)}(t)=\frac{t^n}{n!}e^{-(1+c_q)t}Y(t),
$$ 
as is easily seen by Laplace transformation.

Clearly, $\rho_q\ast\mu$ is a bounded integrable function with
\begin{equation}\label{eq:L}
||\rho_q\ast\mu||_\infty\le ||\rho_q||_\infty||\mu||_1<c_q,\quad
||\rho_q\ast \mu||_1=||\rho_q||_1||\mu||_1<\frac{c_q}{1+c_q},
\end{equation}
and then $\rho_q\ast(\rho_q\ast\mu)^{\ast n},n\ge 1$ is a continuous
integrable function on
$\mathbb R$, vanishing for $t\le n\log(1/q)$ and for
$t\to\infty$. Furthermore,
$$
||\rho_q\ast(\rho_q\ast\mu)^{\ast n}||_\infty < (c_q/(1+c_q))^n,
$$
and this shows that the right-hand side of  \eqref{eq:conv} converges
uniformly on $[0,\infty[$, so $(1-q)\tau_q$ has a continuous density on $[0,\infty[$
tending to 0 at infinity.

For $n\ge 1$ and $x\in [n\log(1/q),\infty[$ we get
\begin{eqnarray*}
\lefteqn{\rho_q^{\ast(n+1)}\ast\mu^{\ast n}(x)}\\ 
&= & \int_0^x
\frac{(x-t)^n}{n!}e^{-(1+c_q)(x-t)}Y(x-t)\,d\mu^{\ast
  n}(t)\\
&=&e^{-(1+c_q)x}\sum_{k=n}^\infty
\frac{(x-k\log(1/q))^n}{n!}q^{-k(1+c_q)}Y(x-k\log(1/q))\mu^{\ast n}(k\log(1/q))
\end{eqnarray*}
which is a finite sum, and
$$
\mu^{\ast n}(k\log(1/q))=\sum_{p_1+\ldots+p_n=k}\prod_{j=1}^n
\frac{q^{2p_j}}{1-q^{p_j}},\;k=n,n+1,\ldots.
$$
In particular, 
$$
\mu^{\ast n}(n\log(1/q))=\left(\frac{q^2}{1-q}\right)^n.
$$
For $n\ge 0$ and $0\le x<(n+1)\log(1/q)$ we then get
\begin{eqnarray}\label{eq:final}
\lefteqn{(1-q)\tau_q(x)=}\\
&&e^{-(1+c_q)x}\sum_{j=0}^n
q^{-j(1+c_q)}Y(x-j\log(1/q))\sum_{k=0}^j\frac{(x-j\log(1/q))^k}{k!}\mu^{\ast
  k}(j\log(1/q)).\nonumber
\end{eqnarray}
 On $[0,\log(1/q)[$ it is equal to
$\exp(-(1+c_q)x)$, on $[\log(1/q),2\log(1/q)[$ it is equal to
$$
\exp(-(1+c_q)x)\left(1+\frac{q^{1-c_q}}{1-q}(x-\log(1/q))\right), 
$$
on $[2\log(1/q),3\log(1/q)[$ it is equal to
\begin{eqnarray*}
\lefteqn{\exp(-(1+c_q)x)\left(1+\frac{q^{1-c_q}}{1-q}(x-\log(1/q))+\right.}\\
&&\left.\frac{q^{2(1-c_q)}}{1-q^2}(x-2\log(1/q))+
\frac{q^{2(1-c_q)}}{2(1-q)^2}(x-2\log(1/q))^2\right),
\end{eqnarray*}
etc.

Using the expression \eqref{eq:final} it is possible to calculate the
derivative of $(1-q)\tau_q$ from the right and from the left at the
point $x=n\log(1/q),n\ge 1$. The difference between the right and the
left derivative equals $q^{2n}/(1-q^n)$ and this gives the jump $J_n$
of \eqref{eq:jump}.
$\square$

It is straightforward to transfer the results of
Theorem~\ref{thm:tauq} to give Theorem~\ref{thm:main} using that
$\nu_q$ is the image measure of $\tau_q$ under $t\mapsto e^{-t}$,
hence $\widetilde{\nu_q}(t)=\tau_q(\log(1/t)),\;0<t\le 1$.

\begin{rem} {\rm The representation \eqref{eq:factor} and
    Theorem~\ref{thm:tauq} show that $(1-q)\tau_q$ is a standard
    $p$-function in the terminology from the theory of regenerative
    phenomena, cf. \cite{K}.}
\end{rem}

\begin{figure}[ht]
\begin{center}
\includegraphics[width=100mm]{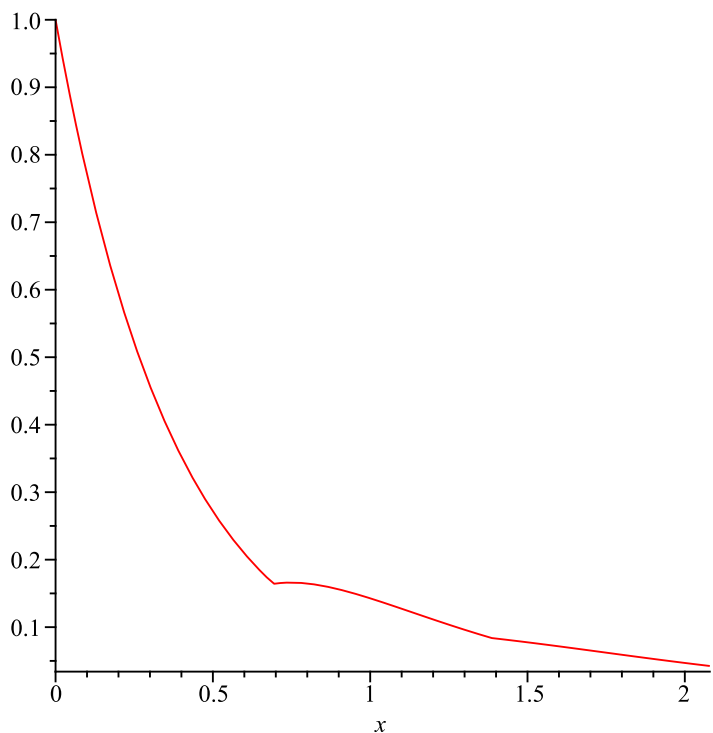}
\caption{The graph of $(1-q)\tau_q$ on $[0,3\log(1/q)]$ for $q=0.5$ }
\label{fig:q5}
\end{center}
\end{figure}

\begin{figure}[ht]
\begin{center}
\includegraphics[width=100mm]{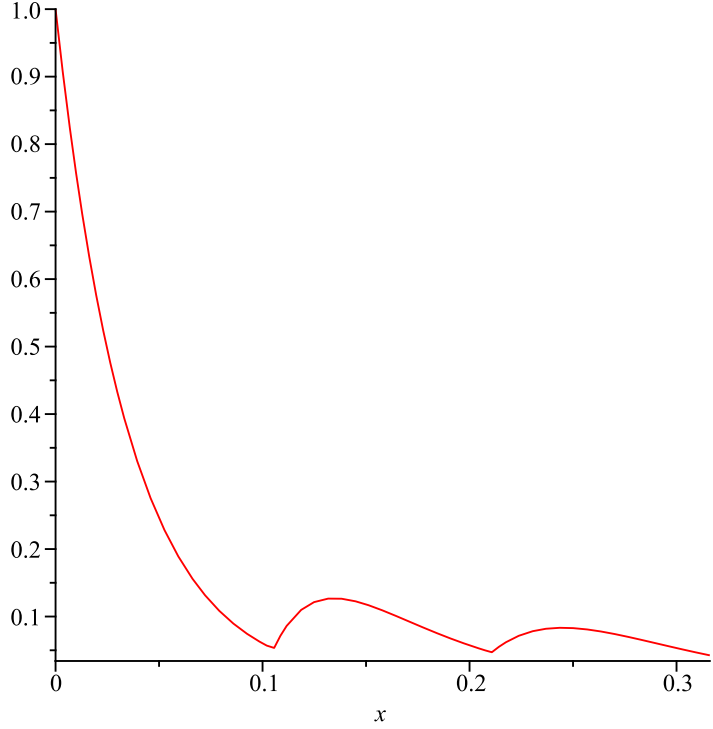}
\caption{The graph of $(1-q)\tau_q$  on $[0,3\log(1/q)]$ for $q=0.9$}
\label{fig:q9}
\end{center}
\end{figure}

\section{Further properties of $\tau_q$}

Formally, by Fourier inversion we get that
$$
\tau_q(x)=\frac{1}{2\pi}\int_{-\infty}^\infty e^{iyx}\frac{dy}{f_q(1+iy)}.
$$
The function $1/f_q(1+iy)$ is a non-integrable $L^2$-function, so the
formula holds in the $L^2$-sense. To see this we notice that
$$
\frac{f_q(z)}{z}=1-q+ \int_0^\infty e^{-tz}h_q(t)\,dt,\quad \Re z>0,
$$
where 
\begin{equation}\label{eq:hq}
h_q(t)=(1-q)\sum_{k>t/\log(1/q)} \frac{q^k}{1-q^k}.
\end{equation} 
In particular
$$
\frac{f_q(1+iy)}{1+iy}=1-q + \int_0^\infty e^{-ity}e^{-t}h_q(t)\,dt,
$$
and since $e^{-t}h_q(t)$ is integrable, it follows from the
Riemann-Lebesgue Lemma that we  get the asymptotic
behaviour
\begin{equation}\label{eq:asymp}
f_q(1+iy) \sim (1-q)(1+iy),\quad |y|\to\infty.
\end{equation}

Furthermore, we notice that
$$
\Re f_q(1+iy)=1-q +(1-q)\sum_{k=1}^\infty\frac{q^k}{1-q^k}\left(1-q^k\cos(ky\log(q))\right),
$$
hence
$$
1\le\Re f_q(1+iy)\le 1-q+\sum_{k=1}^\infty q^k(1+q^k),
$$
showing that $\Re f_q(1+iy)$ is bounded below and above.
It follows that the symmetrized
density

\begin{equation}\label{sym}
\varphi_q(x)=
  \begin{cases}
    \displaystyle \tau_q(x)  & \text{if $x\ge 0$},
\\ &\\
    \displaystyle \tau_q(-x) & \text{if $x<0$},  \end{cases}
\end{equation}
is the Fourier transform of the non-negative integrable function
$$
\frac{2\Re f_q(1+iy)}{|f_q(1+iy)|^2},
$$
and therefore $\varphi_q(x)$ is continuous and positive definite, so
$\tau_q$ is the restriction to $[0,\infty[$ of such a function.

\begin{rem} {\rm The function $f_q$ defined in \eqref{eq:Bdqt} is a
    Bernstein function in the sense of \cite{B:F}, but not a complete
    Bernstein function in the sense of \cite{S:S:V}, because
    $f_q(z)/z$ is not a Stieltjes function as shown by
    formula \eqref{eq:hq}. This is in contrast to 
$$
\lim_{q\to 1}f_q(z)=\psi(z+1)+\gamma,
$$
which is a complete Bernstein function, cf. \cite{B:D1}.
}
\end{rem}

\section{Relation to other work}

The transformation $T$ can be extended from normalized Hausdorff
moment sequences to the set
$\mathcal K=[0,1]^{\mathbb N}$  of
sequences $(x_n)=(x_n)_{n\ge 1}$ of numbers from the unit interval $[0,1]$.
This was done in \cite{B:B}, where  
$T:\mathcal K\to \mathcal K$ is defined by
\begin{equation}\label{eq:trans}
(T(x_n))_n=\frac{1}{1+x_1+\ldots+x_n},\quad n\ge 1.
\end{equation}
The connection is that a normalized Hausdorff moment sequence
$(a_n)_{n\ge 0}$ is considered as the element $(a_n)_{n\ge
  1}\in\mathcal K$.

 Since $T$ is a continuous transformation of the compact convex set
$\mathcal K$
in the space $\mathbb R^{\mathbb N}$ of real sequences
equipped with the product topology, it has a fixed point by
Tychonoff's theorem, and this is $(m_n)_{n\ge 1}$. 

There is no reason a priori that the fixed point $(m_n)$ of
\eqref{eq:trans} should be a
Hausdorff moment sequence, but as we have seen above, the motivation for the study of $T$
comes from the theory of Hausdorff moment sequences.

Although $T$ is not a contraction on $\mathcal K$ in the natural
metric
$$
d((a_n),(b_n))=\sum_{n=1}^\infty 2^{-n}|a_n-b_n|,\quad (a_n),(b_n)\in\mathcal K,
$$
it was proved in \cite{B:B} that $T$ maps $\mathcal K$ into the
compact convex subset 
$$
\mathcal C=\left\{(a_n)\in \mathcal K\mid a_1\ge\tfrac12\right\},
$$
and the restriction of $T$ to $\mathcal C$ is a contraction. It is
therefore possible to infer that $(m_n)$ is an attractive fixed point
from the fixed point theorem of Banach.

\noindent Christian Berg, Helle Bjerg Petersen\\
Department of Mathematical Sciences\\
University of Copenhagen\\
Universitetsparken 5\\
DK-2100 K\o benhavn \O, Denmark\\
E-mail addresses: {\footnotesize  berg@math.ku.dk} (C. Berg)\\
{\footnotesize hellebp@gmail.com} (H. B. Petersen).
 
\end{document}